\documentclass{article}
\usepackage{graphics}
\begin{document}
\sf
\noindent \underline{A system for constructing relatively small polyhedra from Sonob\'e modules.}\footnote{
\sf This article assumes a basic familiarity with the Sonob\'e module, including the 
ability to construct an augmented icosahedron.  [3] (pp. 42ff) and [5] each contain
excellent introductions}

\bigskip

\noindent Adrian Riskin

\noindent Department of Mathematics

\noindent Mary Baldwin College

\noindent Staunton, VA  24401  USA

\noindent ariskin@mbc.edu

\bigskip

Anyone familiar with the Sonob\'e module has probably assembled both the augmented octahedron and the 
augmented icosahedron\footnote{\sf Often erroneously called ``stellated''. A polyhedron is augmented by 
capping each face with a pyramid.} Instructions for both appear in both Kasahara and Takahama and in Mukerji, 
and can also be found on most of the innumerable websites discoverable by googling the word ``Sonob\'e''.
Kasahara and Takahama give some sketchy instructions on how to build some other polyhedra with these modules.
Some of their models are relatively small but their concentration, as does that of the websites, turns 
rapidly to monstrously large polyhedra containing hundreds or even thousands of units.  Thus the casual 
explorer of the wonders of the Sonob\'e module, who perhaps has no desire to fold thousands of units before
construction can even begin, is left without many new models small enough to be built in half an hour or so.

If you have an augmented octahedron at hand, look at it carefully and notice that it consists of eight 
equilateral triangles capped off by low 3-sided pyramids.  The key insight here is that other polyhedra 
all of whose faces are equilateral triangles can also be constructed in their augmented versions out of 
Sonob\'e modules.  Such polyhedra are called \textit{deltahedra} after the Greek capital letter delta: $\Delta$.
The icosahedron is a deltahedron, which is why an augmented version can be constructed of Sonob\'e modules.
It turns out that there are six other convex deltahedra besides the octahedron and the icosahedron [see 1, 6].
These references are useful because they include instructions for making paper models of the deltahedra, which in
turn are useful for designing augmented Sonob\'e versions.  As you become more familiar with the system described
here you will be able to dispense with the paper models.  There are also infinitely many nonconvex deltahedra.
In theory augmented versions of any of them could be constructed from Sonob\'e modules.  In practice, however, 
not all are constructible because steep angles between the faces of a given deltahedron can force the pyramidal
caps to intersect which (unfortunately) isn't possible with standard origami paper.

Now, as you stare at the augmented octahedron, notice that there are two kinds of vertices: The kind at the tips
of the 3-sided pyramids and the kind where four pyramids come together.  There are six of the second kind, 
which correspond to the six original vertices of the octahedron.  For any given deltahedron it is the pattern in
which these original vertices are connected along with the number of pyramids at each one which determine the 
construction of the augmented version.  Given any deltahedron it is possible to draw a diagram which will serve
as a blueprint for the construction of an augmented Sonob\'e version.  The diagram for the octahedron looks like 
this:

\begin{figure}[h]
\begin{center}
\includegraphics{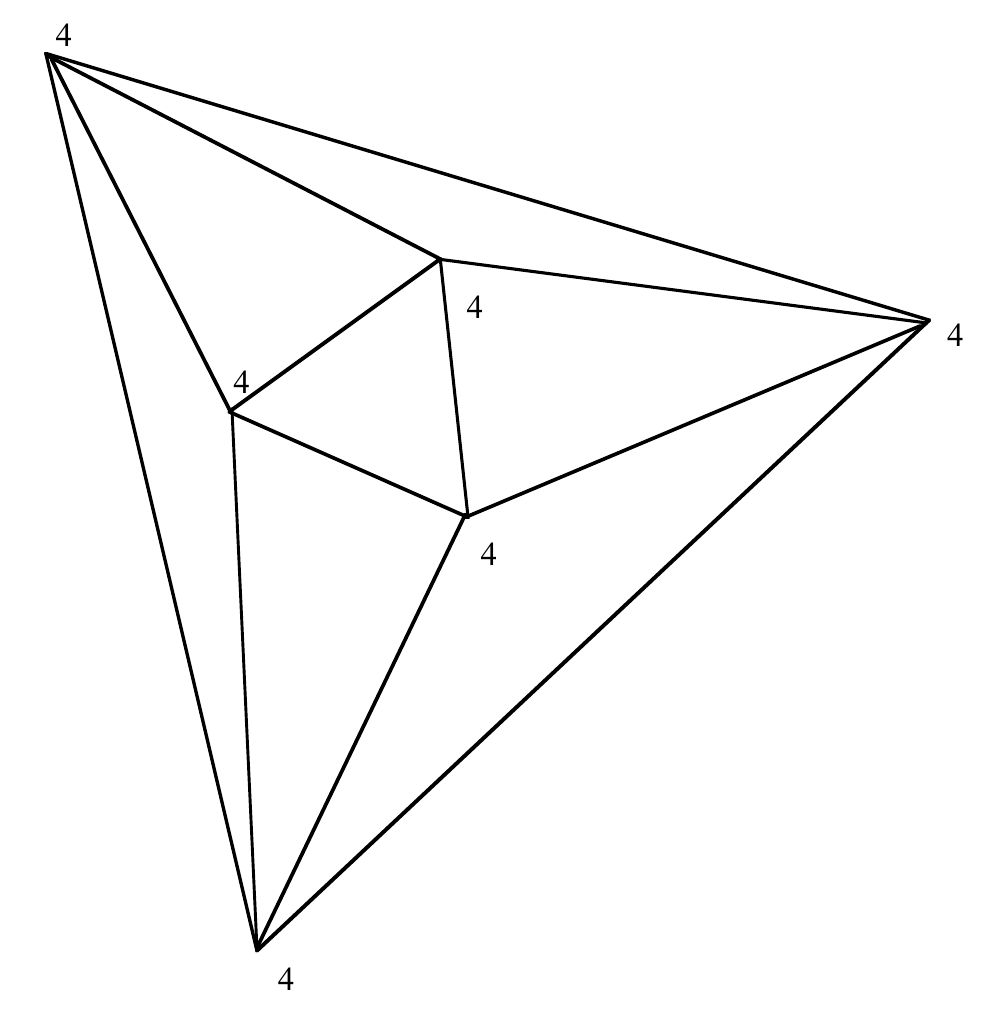}
\end{center}
\end{figure}

Look at the octahedron again.  Notice that each of the original vertices not only has four pyramids meeting at 
it but also has eight edges meeting at it.  Of these eight edges, four are ``mountain'' edges, which go to the 
tip of one of the pyramids, and four are ``valley'' edges, which go to another one of the original vertices in 
the diagram.  One can construct the octahedron from the diagram by choosing an original vertex to start at, 
putting four Sonob\'es in a cycle to represent the four valley edges, and then adding four more modules in the 
usual manner to make four pyramids surrounding the original vertex.  Then choose one of the four original 
vertices connected to the starting vertex by an edge in the diagram, consult the diagram to see that it is
meant to be surrounded by four pyramids, add sufficient Sonob\'es in a cycle around that vertex to make four, 
and then complete the necessary pyramids.  Continue in this manner until the model closes itself off, which will
happen automatically.

Now let's try the method with a new deltahedron: The pentagonal dipyramid.  As far as I can ascertain, instructions
for constructing an augmented version of this polyhedron out of Sonob\'e modules have never appeared either in print 
or on the internet\footnote{\sf Instructions for constructing a modular version of this augmented polyhedron
out of another module are given in [4], but the process is not theorized.}.  Imagine a 5-sided pyramid made of
equilateral triangles with a regular pentagon for a base.  Now imagine sticking two of these together along 
their pentagonal bases to create a deltahedron with ten faces.  It may be possible for you to imagine this clearly
enough to see that the diagram of the original vertices of this deltahedron is:

\begin{figure}[h]
\begin{center}
\includegraphics{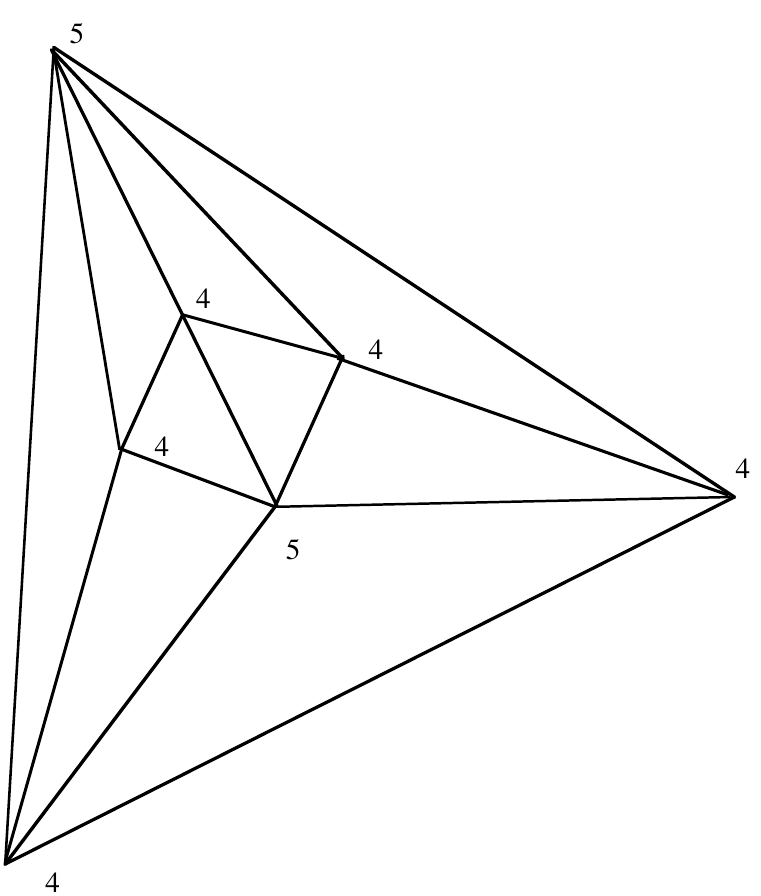}
\end{center}
\end{figure}

Now try constructing the augmented version with modules.  Note that this deltahedron has ten faces.  The number
of modules necessary for constructing an augmented version of a deltahedron is always $\frac{3}{2}$ times the 
number of faces in the original deltahedron, so that fifteen modules will be necessary.  If you find it too
difficult to imagine the polyhedron clearly enough to draw the above diagram it will probably be helpful to make a 
rough paper model of it first.  Enlarge and copy the next diagram, cut it out and fold and tape it into a model
of the pentagonal dipyramid (this kind of diagram is called a \textit{net} of the polyhedron).  Nets for all
the convex deltahedra are available at [6], as well as rotatable pictures of some of the many nonconvex ones.

\begin{figure}[h]
\begin{center}
\includegraphics{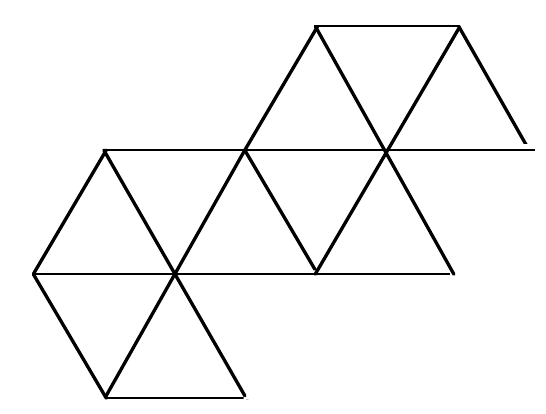}
\end{center}
\end{figure}

By means of this system it is possible to construct models of all the augmented convex deltahedra, as well as
some others.  One particularly interesting nonconvex deltahedron to construct is obtained from a cube by
capping each of the six square faces with a pyramid which has four equilateral triangular sides.  Instructions
for making a beautiful semi-skeletal model of this 24-sided deltahedron can be found in [2].  Another good one
is obtained by augmenting a tetrahedron in the same way --- by capping off each of the four sides with another
tetrahedron.  Finally, any polyhedron whose faces consist of only equilateral triangles and regular hexagons can
be considered to be a deltahedron by dissecting the hexagonal faces into six equilateral triangles.  The simplest
such polyhedron is the truncated tetrahedron, formed from a tetrahedron by cutting off four tetrahedral tips to 
create four equilateral triangles and four regular hexagons.  Augmented polyhedra with hexagonal faces seem to 
be less stable than the others.  That's my method, and here's a challenge: I wonder if anyone can create
augmented toroidal deltahedra from Sonob\'es?  I haven't found a way to do it yet.  Happy folding!

\bigskip

\noindent \underline{References}

\bigskip

\begin{enumerate}
\item Freudenthal, H. and van der Waerden, B. L. "On an Assertion of Euclid." Simon Stevin 25, 115-121, 1947.
\item Fus\'e, Tomoko. Unit Polyhedron Origami.  Japan Publications.  2006.  Tokyo.
\item Kasahara, K. and Takahama, T. Origami for the Connoisseur.  Japan Publications.  1987.  Tokyo.
\item Kasahara, K. Origami Shinhakken 2 (in Japanese).  Nichibou Publishing House 2006.  Tokyo.
\item Mukerji, Meenakshi.  Marvelous Modular Origami.  A. K. Peters.  2007.  Wellesley, Massachusetts, USA.
\item Weisstein, Eric W. "Deltahedron." From MathWorld--A Wolfram Web Resource. 
http://mathworld.wolfram.com/Deltahedron.html
\end{enumerate}

\end{document}